\documentstyle[11pt]{article}
\newlength{\mytopmargin}
\newlength{\myleftmargin}
\newcommand{\al}{\alpha}
\newcommand{\de}{\delta}
\newcommand{\om}{\omega}
\def\rlx{\relax\leavevmode}
\def\inbar{\vrule height1.5ex width.4pt depth0pt}
\def\zz{\rlx\hbox{\small \sf Z\kern-.4em Z}}
\def\rr{\rlx\hbox{\scriptsize \rm I\kern-.18em R}}
\def\nn{\rlx\hbox{\rm I\kern-.18em N}}
\def\qq{\rlx\hbox{\,$\inbar\kern-.3em{\rm Q}$}}
\newcommand{\dif}[1]{\frac{\partial}{\partial #1}}
\newcommand{\bfrac}[2]{\frac{\displaystyle #1}{\displaystyle #2}}

\newcommand{\ph}{\widehat{\Phi}}

\newcommand{\inv}[1]{\widetilde{#1}}
\newtheorem{lemma}{Lemma}[section]

\newtheorem{thm}[lemma]{Theorem}
\newtheorem{cor}[lemma]{Corollary}
\newtheorem{prop}[lemma]{Proposition}
 
\begin{document}
\noindent
\begin{center}{  \Large\bf
A $q$-analogue of the type $A$ Dunkl operator \\[2mm] and integral kernel}
\end{center}
\vspace{5mm}
 
\noindent
\begin{center} T.H.~Baker\footnote{email: tbaker@maths.mu.oz.au; supported
by the ARC} and
  P.J.~Forrester\footnote{email: matpjf@maths.mu.oz.au; supported
  by the ARC}\\[2mm]
{\it Department of Mathematics, University of Melbourne, \\
Parkville, Victoria 3052, Australia}
\end{center}
\vspace{.5cm}
		 
\begin{quote}
We introduce the $q$-analogue of the type $A$ Dunkl operators, which are
a set of degree--lowering operators on the space of polynomials in $n$
variables. This allows the construction of raising/lowering operators with
a simple action on non-symmetric Macdonald polynomials. A bilinear series
of non-symmetric Macdonald polynomials is introduced as a $q$-analogue of 
the type $A$ Dunkl integral kernel ${\cal K}_A(x;y)$. The aforementioned
operators are used to show that the function satisfies $q$-analogues of
the fundamental properties of ${\cal K}_A(x;y)$.
\end{quote}
\date{}

\setcounter{equation}{0}
\section{Introduction}

The purpose of this paper is to obtain $q$-analogues of some fundamental
results concerning type $A$ integral kernels ${\cal K}_A(x;y)$ appearing in 
the works of Dunkl \cite{dunkl89a,dunkl91a,dunkl92a,dunkl96a}. 
The kernel ${\cal K}_A$ allows a multidimensional analogue of the Fourier
transform to be constructed, and plays a pivotal role in recent studies
of generalized Hermite polynomials \cite{forr96c,forr96d}. A necessary
prerequisite will be to give a suitable definition of the $q$-analogue 
of the type $A$ Dunkl operator
\begin{equation} \label{dunklop}
d_i := \frac{\partial}
{\partial x_i} + \frac{1}{\alpha}\sum_{j\neq i} \frac{1}{x_i-x_j}
(1-s_{ij})
\end{equation}
where $s_{ij}$ acts on functions of $x:=(x_1,\ldots, x_n)$ by interchanging 
the variables $x_i$ and $x_j$.

For the purpose of comparison with later results, we recall some 
results concerning the kernel ${\cal K}_A(x;y)$ and the Dunkl operator $d_i$. 
The former is a bilinear series 
in non-symmetric Jack polynomials, denoted $E_{\eta}(x)$, which themselves
are eigenfunctions of the Cherednik operators \cite{cher91a}
\begin{equation} \label{cherednik.1}
\xi_i = \al x_i\dif{x_i} + \sum_{p<i} \frac{x_i}{x_i-x_p}(1-s_{ip})
+ \sum_{p>i} \frac{x_p}{x_i-x_p}(1-s_{ip}) +1-i
\end{equation}
The Dunkl operator $d_i$ is related to the Cherednik operator
$\xi_i$ via
\begin{equation} \label{duncher}
d_i = \frac{1}{\alpha x_i}\left( \xi_i +n-1 - \sum_{p>i} s_{ip}
\right)
\end{equation}
and this will be our starting point in defining an suitable
$q$-analogue of $d_i$.

Following Sahi \cite{sahi96a}, for a node $s=(i,j)$ in a composition
$\eta:=(\eta_1,\eta_2,\ldots,\eta_n)\in\nn^n$, define 
the arm length $a(s)$, arm colength $a'(s)$, leg length
$l(s)$ and  leg colength $l'(s)$ by
\begin{eqnarray}
a(s)= \eta_i - j && l(s) = \#\{k>i|j\leq \eta_k\leq\eta_i\} \;+\;
\#\{k<i|j\leq \eta_k+1\leq\eta_i\} \nonumber\\
a'(s)=j - 1 && l'(s) = \#\{k>i| \eta_k > \eta_i\} \;+\;
\#\{k<i|\eta_k\geq\eta_i\}  \label{guion}
\end{eqnarray}
Using these, define constants 
\begin{eqnarray}
d_{\eta} &:=& \prod_{s\in\eta} (\al(a(s)+1) + l(s)+1)\hspace{2cm} 
d'_{\eta} := \prod_{s\in\eta} (\al(a(s)+1) + l(s)) \nonumber\\
e_{\eta} &:=& \prod_{s\in\eta} (\al(a'(s)+1) + n-l'(s)) \label{cajamarca}
\end{eqnarray}
With these constants, the type $A$ kernel
is defined as \cite{forr96d}
\begin{equation} \label{ker-A}
{\cal K}_A(x;y) = \sum_{\eta} \al^{|\eta|}\;\frac{d_{\eta}}
{d'_{\eta}e_{\eta}}\;E_{\eta}(x)\:E_{\eta}(y)
\end{equation}

Set $s_i:=s_{i,i+1}$ for $1\leq i\leq n-1$.  The following 
raising/lowering operators 
\begin{eqnarray}
\Phi &:=& x_n\, s_{n-1}\cdots s_2\, s_1 = s_{n-1}\cdots s_i \:x_i \:
s_{i-1}\cdots s_1 \label{jphi}\\
\ph &:=& d_1 s_1 s_2 \cdots s_{n-1} = s_1s_2\cdots s_{i-1}\; d_i\;
s_i s_{i+1}\cdots s_{n-1} \label{jph}
\end{eqnarray}
have a very simple action on the non-symmetric Jack polynomials
$E_{\eta}(x)$ \cite{knop96c,forr96d},
\begin{eqnarray}
\Phi\, E_{\eta} &=& E_{\Phi\eta} \label{ica.1}\\
\ph\, E_{\eta} &=& \frac{1}{\alpha} \frac{d'_{\eta}}{d'_{\ph\eta}}\,
E_{\ph\eta} \label{ica.2}
\end{eqnarray}
where $\Phi\eta := (\eta_2,\eta_3,\ldots,\eta_n,\eta_1+1)$ and
$\ph\eta := (\eta_n -1,\eta_1,\eta_2,\ldots,\eta_{n-1})$.

The fundamental properties of the kernel ${\cal K}_A(x;y)$ are given 
by the following result \cite[Theorem 3.8]{forr96d}

\begin{thm} \label{jackmain}
The function ${\cal K}_A(x;y)$ possesses the following properties
\begin{eqnarray*}
(a)\qquad s_i^{(y)}\;{\cal K}_A(x;y) &=& s_i^{(x)}\;{\cal K}_A(x;y) \\
(b)\qquad \ph^{(y)}\;{\cal K}_A(x;y) &=& \Phi^{(x)}\;{\cal K}_A(x;y) \\
(c)\qquad d_i^{(y)}\;{\cal K}_A(x;y) &=& x_i\;{\cal K}_A(x;y)
\end{eqnarray*}
\end{thm}

The above kernel has a symmetric counterpart ${}_0{\cal F}_0(x;y)$
which itself is expressed in terms of the symmetric Jack polynomials
$P^{(\alpha)}_{\lambda}(x)$ \cite{mac}. The symmetric Jack polynomials
can be expressed in terms of their non-symmetric siblings $E_{\eta}(x)$
by \cite{sahi96a}
\begin{equation}
P^{(\alpha)}_{\kappa}(x) = d'_{\kappa}\sum_{\eta }
\frac{1}{d'_{\eta}}\,E_{\eta}(x) \label{kirsten.1}
\end{equation}
where the sum is over distinct permutations $\eta$ of the partition 
$\kappa$. They can also be obtained by symmetrization \cite{forr96d}
\begin{equation}\label{kirsten.2}
{\rm Sym}\: E_{\eta}(x) = 
\frac{n!\, e_{\eta}}{d_{\eta}P^{(\alpha)}_{\lambda}(1^n)}
P^{(\alpha)}_{\eta^+}(x)
\end{equation}
where ${\rm Sym}$ denotes the operation of symmetrization of the variables
$x_1,\ldots,x_n$ and $\eta^+$ denotes the (unique) partition associated
with $\eta$, obtained by permuting its entries. It was shown in \cite{forr96d}
that the symmetric kernel
$$
{}_0{\cal F}_0(x;y) := \sum_{\kappa} \frac{\alpha^{|\kappa|}}{d'_{\kappa}}
\frac{P^{(\alpha)}_{\kappa}(x)\,P^{(\alpha)}_{\kappa}(y)}
{P^{(\alpha)}_{\kappa}(1^n)}
$$
can be obtained from the non-symmetric kernel ${\cal K}_A(x;y)$ via 
symmetrization:
\begin{equation} \label{belem}
{\rm Sym}^{(x)}\:{\cal K}_A(x;y) = n!\,{}_0{\cal F}_0(x;y).
\end{equation}

In this work, we shall be concerned with providing $q$-analogues of the
above results. After introducing preliminary results and notations 
dealing with non-symmetric Macdonald polynomials $E_{\eta}(x;q,t)$
and (type $A$) affine Hecke algebras, we proceed to define an analogue
of the Dunkl operator (\ref{dunklop}) and show that they form a mutually
commuting set of degree-lowering operators. We then
construct the analogue of Knop and Sahi's
raising operator $\Phi$ given by (\ref{jphi}), as well as its lowering
counterpart $\ph$, demonstrating their simple action on $E_{\eta}(x;q,t)$. 
Finally, we construct a $q$-analogue of the kernel ${\cal K}_A(x;y)$ and 
derive the corresponding version
of Theorem \ref{jackmain}. Symmetrization of this kernel is then shown to
recover the well-known symmetric version ${}_0{\cal F}_0(x;y;q,t)$
\cite{macunp1,lass97a}.

\setcounter{equation}{0}
\section{Preliminaries}

We begin by presenting the standard realization of the affine Hecke
algebra on the space of polynomials in $n$ variables (see e.g. 
\cite{noumi96a,vinet96b}).
Let $\tau_i$ be the $q$-shift operator in the variable $x_i$, so
that 
$$
(\tau_i\,f)(x_1,\ldots,x_i,\ldots,x_n) = 
f(x_1,\ldots,qx_i,\ldots,x_n).
$$
The Demazure-Lustig operators are defined by
\begin{eqnarray} \label{pucallpa}
T_i &=& t + \frac{tx_i-x_{i+1}}{x_i - x_{i+1}}\left(s_i -1 \right)
\hspace{2cm} i=1,\ldots,n-1 \\
\mbox{and} \hspace{2cm}
T_0 &=& t + \frac{qtx_n-x_1}{qx_n - x_1}\left(s_0 -1 \right)
\end{eqnarray}
where $s_0 := s_{1n}\tau_1\tau_n^{-1}$.
For future reference, we note the following action of $T_i$, 
$1\leq i\leq n-1$ on the monomial $x_i^a\,x_{i+1}^b$
\begin{equation} \label{t.action}
T_i\, x_i^a x_{i+1}^b = \left\{ \begin{array}{ll}
(1-t)x_i^{a-1}x_{i+1}^{b+1} + 
\cdots +(1-t)x_i^{b+1}x_{i+1}^{a-1} + x_i^bx_{i+1}^a& a > b \\
t x_i^ax_{i+1}^a & a=b \\
(t-1)x_i^{a}x_{i+1}^{b} + 
\cdots +(t-1)x_i^{b-1}x_{i+1}^{a+1} + t x_i^bx_{i+1}^a & a < b 
\end{array} \right.
\end{equation}

In addition to the operators $T_i$, define
$$
\om := s_{n-1}\cdots s_2\,s_1\tau_1 = 
s_{n-1}\cdots s_i\tau_i s_{i-1}\cdots s_1 .
$$
The affine Hecke algebra is then generated by elements $T_i$, 
$0\leq i\leq n-1$ and $\om$, satisfying the relations
\begin{eqnarray}
(T_i-t)\,(T_i+1) &=& 0 \label{tdefs.1}\\
T_i\;T_{i+1}\;T_i &=& T_{i+1}\;T_i\;T_{i+1} \label{tdefs.2}\\
T_i\;T_j &=& T_j\;T_i \qquad |i-j| \geq 2 \\
\om\;T_i &=& T_{i-1}\;\om  \label{tdefs.4}
\end{eqnarray}

{}From the quadratic relation (\ref{tdefs.1}), we have the identity
\begin{equation} \label{t-inv}
T_i^{-1} = t^{-1} - 1 +t^{-1}\,T_i .
\end{equation}
Useful relations between the operators $\om$, $T_i$ and $x_i$,
$1\leq i\leq n-1$, include
\begin{eqnarray}
T_i^{-1}\,x_{i+1} =t^{-1}x_i\,T_i \hspace*{2cm} &
T_i^{-1}\,x_i = x_{i+1}\,T_i^{-1} + (t^{-1}-1)x_i  \label{id.1}\\
T_i\,x_i = t x_{i+1}\,T_i^{-1} \hspace*{2cm} & 
T_i\,x_{i+1} = x_i\,T_i +(t-1)x_{i+1} \label{id.2}\\
\om\,x_i = qx_n\om \hspace{2cm} & \om\,x_{i+1} = x_i\,\om \quad
1\leq i\leq n-1 \label{id.3}
\end{eqnarray}

We can define operators\footnote{The normalization is chosen such that
as $q\rightarrow 1$, $(1-Y_i)/(1-q) \rightarrow \xi_i/\alpha$ with
$\xi_i$ given by (\ref{cherednik.1})} \cite{noumi96a}
\begin{equation}
Y_i = t^{-n+i}\;T_i\cdots T_{n-1}\;\om\;T_1^{-1}\cdots T_{i-1}^{-1},
\hspace{2cm} 1\leq i\leq n
\end{equation}
which commute amongst themselves: $[Y_i,Y_j]=0$, $1\leq i,j\leq n$.
They also possess the following relations with the operators $T_i$
\begin{equation}\label{yt1}
T_i\,Y_{i+1}\,T_i = t\,Y_i \hspace{2cm} [T_i, Y_j] =0, \quad j\neq i, i+1
\end{equation}
while the following relations with $x_n$ will be needed in Section 4
\begin{eqnarray}
Y_i\,x_n &=& x_n\,Y_i + t^{-n+1}(1-t)x_n\,T_i\cdots T_{n-1}\omega
T_1^{-1}\cdots T_{i-1}^{-1} \qquad 1\leq i\leq n-1 \label{llave.1}\\
Y_n\,x_n &=& qt^{-n+1} x_n\omega T_1\cdots T_{n-1} \label{llave.2}
\end{eqnarray}
These identities follow from a direct calculation involving
(\ref{id.3}), (\ref{id.2}) and (\ref{t-inv}).

The fact that the operators $Y_i$ mutually commute implies that they possess
a set of simultaneous eigenfunctions, the non-symmetric Macdonald
polynomials. Specifically, let $\prec$ denote the partial order on
compositions $\eta$ defined for $\eta\neq\nu$ by
$$
\nu\prec\eta \quad\mbox{iff}\quad \nu^+ < \eta^+
\quad\mbox{or in the case $\nu^+ = \eta^+$}\quad \nu < \eta
$$
where $<$ is the usual dominance order for
$n$-tuples, i.e. $\nu < \eta$ iff $\sum_{i=1}^p (\eta_i - \nu_i) \geq 0$,
for all $1\leq p \leq n$.
Then the non-symmetric Macdonald polynomials $E_{\eta}(x;q,t)$ can 
be defined by the conditions 
\begin{eqnarray*}
E_{\eta}(x;q,t) &=& x^{\eta} + \sum_{\nu\prec \eta} b_{\eta\nu} x^{\nu} \\
Y_i \, E_{\eta}(x;q,t) &=& t^{\bar{\eta}_i}\,E_{\eta}(x;q,t) \qquad
1\leq i \leq n
\end{eqnarray*}
where 
\begin{equation}\label{e-val.1}
\bar{\eta}_i = \al\eta_i - \#\{k<i\;|\;\eta_k \geq \eta_i\} -
\#\{k>i\;|\;\eta_k > \eta_i\}
\end{equation}
with $\alpha$ a parameter such that $t^{\alpha}=q$.

Define the analogue of the constants (\ref{cajamarca}) by \cite{noumi96c}
\begin{eqnarray}
d_{\eta}(q,t) &:=& \prod_{s\in\eta} \left( 1-q^{a(s)+1}t^{l(s)+1} \right)
\quad
d'_{\eta}(q,t):=\prod_{s\in\eta} \left( 1-q^{a(s)+1}t^{l(s)} \right)
\quad \nonumber\\
e_{\eta}(q,t) &:=& \prod_{s\in\eta} \left( 1-q^{a'(s)+1}t^{n-l'(s)} \right)
\label{constants}
\end{eqnarray}
Certain properties of these coefficients 
follow immediately from \cite[Lemmas 4.1, 4.2]{sahi96a} 

\begin{lemma} \label{props}
We have
\begin{eqnarray*}
\frac{d_{\Phi\eta}(q,t)}{d_{\eta}(q,t)} = 
\frac{e_{\Phi\eta}(q,t)}{e_{\eta}(q,t)} = 1-qt^{n+\bar{\eta}_1},
\quad
\frac{d'_{\Phi\eta}(q,t)}{d'_{\eta}(q,t)} = 1-qt^{n-1+\bar{\eta}_1}, \quad
e_{s_i\eta}(q,t) = e_{\eta}(q,t), \\
\frac{d_{s_i\eta}(q,t)}{d_{\eta}(q,t)} = \frac{1-t^{\de_{i,\eta}+1}}
{1-t^{\de_{i,\eta}}}, \qquad
\frac{d'_{s_i\eta}(q,t)}{d'_{\eta}(q,t)} = \frac{1-t^{\de_{i,\eta}}}
{1-t^{\de_{i,\eta}-1}} \qquad \mbox{{\rm for} $\eta_i > \eta_{i+1}$}
\end{eqnarray*}
\end{lemma}

\setcounter{equation}{0}
\section{$q$-Dunkl operators}

We introduce the $q$-Dunkl operators $D_i$ for $1\leq i\leq n$ according to
\begin{equation}\label{qdunkl-def}
D_i := x_i^{-1}\;\left( 1-t^{n-1}\left[ 1 +(t^{-1}-1)\sum_{j=i+1}^n
t^{j-i}\,T_{ij}^{-1} \right] Y_i \right)
\end{equation}
where for $i<j$,
\begin{eqnarray*}
T_{ij}^{-1} &:=& T_i^{-1}\;T_{i+1}^{-1}\;\cdots\;T_{j-2}^{-1}\;T_{j-1}^{-1}
\;T_{j-2}^{-1}\; \cdots\; T_i^{-1}   \\
&=& T_{j-1}^{-1}\; T_{j-2}^{-1}\;\cdots\;T_{i+1}^{-1}\;T_i^{-1}\;
T_{i+1}^{-1}\;\cdots\;T_{j-1}^{-1}
\end{eqnarray*}
Note that as $q\rightarrow 1$ in (\ref{qdunkl-def}), since $T_{ij}^{-1}
\rightarrow s_{ij}$ and $(1-Y_i)/(1-q)\rightarrow \xi_i/\alpha$, we recover 
the type $A$ Dunkl operators due to the relation (\ref{duncher}):
$$
\lim_{q\rightarrow 1} \frac{D_i}{1-q} = d_i
$$

The relations between the operators $D_i$ and the elements $T_i$, $\omega$
of the affine Hecke algebra are given by the following two lemmas,

\begin{lemma} \label{lemma.1}
\begin{eqnarray}
T_i\,D_{i+1} &=& t\,D_i\,T_i^{-1}, \qquad T_i\,D_i = D_{i+1}\,T_i 
+ (t-1)D_i \hspace{2cm} 1\leq i \leq n-1 \label{sapo.1}\\
{}[T_i, D_j] &=& 0 \hspace{9cm} j\neq i, i+1 \label{sapo.2}
\end{eqnarray}
\end{lemma}
{\it Proof.}\quad First note that the second relation in (\ref{sapo.1})
follows from the first relation by multiplying the latter on the left
by $T_i^{-1}$, on the right by $T_i$, and then using (\ref{t-inv}).

{}From (\ref{t-inv}) it follows that $T_i^{-1}$ obeys
the quadratic relation
\begin{equation} \label{qti}
T_i^{-2} + (1 - t^{-1}) T_i^{-1} -t^{-1} = 0 .
\end{equation}
Multiply the relation $t^{-1}Y_{i+1} = T_i^{-1} Y_i T_i^{-1}$ (which 
follows directly from (\ref{yt1})) on the left by $T_i^{-1}$ and apply
(\ref{qti}) to give
\begin{equation}\label{yt2}
T_i^{-1}\, Y_{i+1} = (t^{-1}-1)\,Y_{i+1} + Y_i\,T_i^{-1} .
\end{equation}
{}From the definition (\ref{qdunkl-def}) and the relations 
$T_i x_{i+1}^{-1} = t x_i^{-1} T_i^{-1}$ (which follows from the first
equation in (\ref{id.2})), (\ref{yt1}) and (\ref{yt2}), the first
relation in (\ref{sapo.1}) can be deduced.

Turning to (\ref{sapo.2}), note that a convenient representation of
$D_j$ for $j<n$ in terms of 
\begin{equation} \label{explicit}
D_n:=x_n^{-1}(1-t^{n-1}Y_n)
\end{equation}
is
\begin{equation} \label{crep}
D_j = t^{-n+j}\:T_j\,\cdots\, T_{n-1}\,D_n\,T_{n-1}\,\cdots\,T_j .
\end{equation}
Thus for $i<j-1$ it follows that $[T_i, D_j]=0$. For $j-1 \leq i \leq n-1$
we have
\begin{eqnarray*}
T_i\,D_j &=& t^{-n+j} T_j\cdots T_{n-1}T_{i-1}\,D_n\,T_{n-1}\cdots T_j \\
&=& t^{-n+j} T_j\cdots T_{n-1}\,D_n\,T_{i-1},T_{n-1}\cdots T_j \\
&=& D_j\,T_i , 
\end{eqnarray*}
where the first equality follows from (\ref{tdefs.2}), the second from the 
already established commutativity of $T_i$, $D_j$ for $i<j-1$, and the
final equality follows by further use of (\ref{tdefs.1}).  \hfill $\Box$

\begin{lemma} \label{lemma.2}
We have
\begin{eqnarray}
\omega\,D_{i+1} &=& D_i\,\omega \hspace{4cm} 1\leq i \leq n-1 
\label{horma.1}\\
q\omega\,D_1 &=& D_n\,\omega \label{horma.2}
\end{eqnarray}
\end{lemma}
{\it Proof.}\quad We first prove (\ref{horma.1}) in the special case 
$i=n-1$. To this end, note that for $i\geq 2$,
$$
\omega\,Y_i = Y_{i-1}\,\omega\;+\;(1-t) T_{i-1}\cdots T_{n-2} T_{n-1}^{-1}
T_{n-2}^{-1} \cdots  T_{i-1}^{-1}\,Y_{i-1}\,\omega
$$
which follows from using (\ref{tdefs.4}) to shift the operator $\omega$
to the right. In particular
\begin{equation} \label{borrar}
\omega\,Y_n = Y_{n-1}\,\omega\;+\;(1-t)T_{n-1}^{-1}\,Y_{n-1}\,\omega  .
\end{equation}
The use of (\ref{id.3}) and (\ref{borrar}) and the explicit expression
for $D_n$ given by (\ref{explicit}) allows one to show that
$\omega\,D_n = D_{n-1}\,\omega$. For the cases $i<n-1$ the result 
follows from the case $i=n-1$ and the representation of $D_i$ in terms 
of $D_n$ given by (\ref{crep}) since
\begin{eqnarray*}
\omega\,D_{i+1} &=& t^{-n+i+1}\omega T_{i+1}\cdots T_{n-1}D_n T_{n-1}
\cdots T_{i+1} \\
&=& t^{-n+i+1} T_i \cdots T_{n-2} D_{n-1} T_{n-2} \cdots T_i\:\omega \\
&=& t^{-n+i+1} T_i \cdots T_{n-2}\left(t^{-1}T_{n-1}D_n T_{n-1} \right)
T_{n-2} \cdots T_i\,\omega  = D_i\,\omega
\end{eqnarray*}

To prove (\ref{horma.2}), first note that repeated use of (\ref{qti})
yields that, for $i<j$
\begin{eqnarray}
I_{ij}^{-1} &:=& T_i^{-1}\,T_{i+1}^{-1}\,\cdots\, T_j^{-1}\,T_j^{-1}\,
\cdots\, T_{i+1}^{-1}\,T_i^{-1} \nonumber \\ \label{huelga}
&=& t^{i-j-1} + (t^{-1}-1)\,\sum_{p=i+1}^{j+1} t^{p-j-1}\,T_{ip}^{-1}
\end{eqnarray}
It then follows that 
$$
Y_n\,\omega = \omega\,T_1^{-1}\cdots T_{n-1}^{-1}\,\omega 
= t^{n-1}\omega\,I_{1,n-1}^{-1}\,Y_1 
$$
and so from this and (\ref{huelga}), we have
\begin{eqnarray*}
D_n\,\omega &=& x_n^{-1} \left(1-t^{n-1}Y_n \right)\,\omega
= x_n^{-1} \omega\left(1-t^{2n-2} I_{1,n-1}^{-1} Y_1 \right) \\
&=& q\omega\,x_1^{-1} \left( 1 -t^{n-1} \left[ 1 +(t^{-1}-1)
\sum_{p=2}^n t^{p-1}\,T_{1p}^{-1} \right] Y_1 \right) 
= q\omega\,D_1
\end{eqnarray*}
\hfill $\Box$
\vspace{3mm}\\
{\it Remarks.}\\
{\it 1.}\quad The final relations between the operators $D_i$, 
$1\leq i\leq n$ and the generators of the affine Hecke algebra are the
ones involving the generator $T_0$. These takes the form
\begin{eqnarray*}
T_0\,D_1 &=& q^{-1}t D_n\,T_0^{-1}, \hspace{2cm}
T_0\,D_n = q D_1\,T_0 + (t-1)D_n \\
\left[T_0, D_i\right] &=& 0, \hspace{4cm} 2\leq i\leq n-1
\end{eqnarray*}
which follow immediately using the fact that $T_0 = \om\,T_1\,\om^{-1}$
along with Lemmas \ref{lemma.1} and \ref{lemma.2}.
\vspace{1mm}

\noindent
{\it 2.}\quad It follows from (\ref{huelga}) that 
\begin{eqnarray}
D_i &=& x_i^{-1}\left( 1 -t^{2n-i-1}\,I_{i,n-1}^{-1}\,Y_i\right)
\nonumber\\
&=& x_i^{-1} \left( 1-t^{n-1}\,T_i^{-1}\cdots T_{n-1}^{-1}\,\omega\,
T_1^{-1}\cdots T_{i-1}^{-1} \right)
\end{eqnarray}
providing an alternative definition of the $q$-Dunkl operators.
\vspace{1mm}

Another set of relations which shall be needed later on is
an analogue of \cite[Lemma 3.1]{forr96c}

\begin{lemma} \label{caliente}
We have
\begin{eqnarray}
[D_i, Y_j] &=& \left\{\begin{array}{cc}
t^{i-j}(1-t) Y_j\,T_{ij}\,D_j \qquad i < j \\[2mm]
t^{j-i}(1-t) Y_i\,T_{ji}\,D_i \qquad i > j  \end{array} \right. 
\label{piel.1}\\
D_i\,Y_i - qY_i\,D_i &=& (t-1)\sum_{p=i+1}^n t^{-p+i} Y_p\,T_{ip}\,D_p
+ q(t-1)\,Y_i\, \sum_{p=1}^{i-1}t^{p-i}\,T_{ip}\, D_i \label{piel.2}
\end{eqnarray}
\end{lemma}
{\it Proof.}\quad We start with (\ref{piel.1}) when $i<j$. In this case,
we can use Lemmas \ref{lemma.1} and \ref{lemma.2} to shuffle the $D_i$
to the right to get 
\begin{eqnarray*}
D_i\,Y_j &=& t^{-n+j} T_j\cdots T_{n-1}\,\omega\,T_1^{-1}\cdots T_{i-1}^{-1}
\,D_{i+1}\,T_i^{-1} T_{i+1}^{-1}\cdots T_{j-1}^{-1}  \\
&=& Y_j\,D_i\;+\;t^{-n+j}(t^{-1}-1)\, T_j\cdots T_{n-1}\,\omega\,
T_1^{-1}\cdots T_{i-1}^{-1}\,D_{i+1}\,T_{i+1}^{-1}\cdots T_{j-1}^{-1}
\end{eqnarray*}
where we have used the fact that $D_{i+1}T_i^{-1} = T_i^{-1}D_i + (t^{-1}-1)
D_{i+1}$ (which follows directly from (\ref{sapo.1}) ). Using Lemma
\ref{lemma.1} on the second of these terms to move $D_{i+1}$ to the right
results in the term $t^{i-j}(1-t) Y_j\,T_{ij}\,D_j$ as required. The proof
in the case $i>j$ is somewhat similar.

To prove (\ref{piel.2}) we must also consider 2 cases: $1\leq i<n$ and $i=n$.
For the case $1\leq i<n$, we use the identity
$$
D_i\,T_i = T_i\,D_{i+1} + (t-1)D_i
$$
(which can be obtained from Lemma \ref{lemma.1} ) to move the operator
$D_i$ to the right in the expression
$$
D_i\,Y_i = t^{-n+i} D_i\,T_i\cdots T_{n-1}\,\omega\,T_1^{-1} \cdots
T_{i-1}^{-1} 
$$
and thus obtain
$$
D_i\,Y_i = (t-1)\sum_{p=i+1}^n t^{-p+i} Y_p\,T_{ip}\,D_p +
qt^{-i+1} Y_i\:T_{i-1}\cdots T_1 T_1 \cdots T_{i-1}\:D_i
$$
The second term in the above expression can be simplified using a
result similar to (\ref{huelga}), namely that for $i<j$
$$
T_j\,T_{j-1}\,\cdots\,T_i\,T_i\,\cdots T_{j-1}T_j = 
t^{j-i+1} \;+\; (t-1)\sum_{p=i}^j t^{p-i} T_{p,j+1}
$$
and the stated result follows. The case $i=n$ is derived similarly.
\hfill $\Box$
\vspace{3mm}

In the rest of this section we shall show that the $q$-Dunkl operators $D_i$
commute amongst themselves. We do this by showing that all the $D_i$
commute with $D_n$ (recall that $D_n$ has the simplest form amongst 
all the $q$-Dunkl operators) from which the general result follows
swiftly.

First note that for $i<n$
\begin{equation} \label{agua}
Y_i\,x_n^{-1} = x_n^{-1}\,Y_i + t^{n-i-1}(t-1)x_i^{-1}\,
T_{in}^{-1}\,Y_i
\end{equation}
which follows from pulling $x_n^{-1}$ to the left using 
\begin{equation} \label{suelo}
T_i\,x_{i+1}^{-1} = t x_i^{-1}\,T_i^{-1} \hspace{2cm}
T_i\,x_i^{-1} = x_{i+1}^{-1}\,T_i + (t-1)x_i
\end{equation}
(which themselves follow from (\ref{id.1}) and (\ref{id.2}) ).
Using (\ref{agua}) and (\ref{suelo}), a series of manipulations
yields the relation
$$
Y_n\,T_{in}\,x_n^{-1} = t^{2(n-i)-1}\,x_i^{-1}\,T_{in}^{-1}\,Y_i
$$
Using this and Lemma \ref{caliente} we can rewrite the commutator $[D_i,Y_n]$
in the form
\begin{equation} \label{montreal}
\left[D_i,Y_n \right] = t^{n-i}(t^{-1}-1)\,x_i^{-1}\,T_{in}^{-1}\,Y_i\,
(1-t^{n-1})\,Y_n
\end{equation}

Another result we need is
\begin{eqnarray}
\left[T_{in}^{-1}, x_n^{-1} \right] &=& (x_i^{-1}-x_n^{-1})\,T_{in}^{-1}
+ (t^{-1}-1)\,x_n^{-1}\, I_{i,n-1}^{-1} \nonumber \\
&&+\;(t^{-1}-1)\sum_{p=i+1}^{n-1}x_p^{-1}\,T_i^{-1}\cdots T_{p-2}^{-1} 
T_p^{-1} \cdots T_{n-1}^{-1} \cdots T_i^{-1} \label{reloj}
\end{eqnarray}
which can be derived by repeated use of 
$$
T_i^{-1}x_{i+1}^{-1} = x_i^{-1}\,T_i^{-1} + (t^{-1}-1)x_{i+1}^{-1}
$$
(which itself is derived from (\ref{suelo}) ). 

This is used in the derivation of the final necessary ingredient
\begin{lemma} \label{frio}
We have for $i<n$
$$
\left[ D_i, x_n^{-1} \right] = t^{2n-i-1}(t^{-1}-1)\, x_i^{-1}\,x_n^{-1}
\, T_{in}^{-1}\, Y_i
$$
\end{lemma}
{\it Proof.}\quad  Write $D_i = x_i^{-1}\,(1-t^{n-1}A_i\,Y_i)$ where
$$
A_i := 1 + (t^{-1}-1) \sum_{j=i+1}^n t^{j-i}\,T_{ij}^{-1}
$$
Then we have
\begin{eqnarray}
\left[ D_i, x_n^{-1} \right] &=& \left[ x_i^{-1}\,(1-t^{n-1}A_i\,Y_i),
x_n^{-1} \right] \nonumber \\
&=& -t^{n-1}x_i^{-1} \left( \left[A_i, x_n^{-1} \right] Y_i
+ A_i \left[ Y_i, x_n^{-1} \right] \right) \label{boli}
\end{eqnarray}
Note that
\begin{equation} \label{gafas}
\left[A_i, x_n^{-1} \right] = (t^{-1}-1)t^{n-i} \left[ T_{in}^{-1},
x_n^{-1} \right] 
\end{equation}
while from (\ref{agua}) we have
\begin{equation} \label{lapiz}
A_i\left[ Y_i, x_n^{-1} \right] = (1-t^{-1}) t^{n-i} \left(
x_i^{-1}\,T_{in}^{-1} + (t^{-1}-1) \sum_{j=i+1}^n t^{j-i}\,
T_{ij}^{-1}\,x_i^{-1}\,T_{in}^{-1} \right) Y_i
\end{equation}
However, for $i<j$ application of the relation $T_i^{-1}x_i^{-1} =
t^{-1}x_{i+1}^{-1} T_i$ tells us that
$$
T_{ij}^{-1}\, x_i^{-1} = t^{-j+i} x_j^{-1}\, T_i^{-1}\cdots T_{j-2}^{-1}
T_{j-1} T_{j-2}\cdots T_i
$$
so that 
$$
T_{ij}^{-1}\,x_i^{-1}\,T_{in}^{-1} = \left\{ \begin{array}{ll}
t^{-j+i} x_j^{-1}\,T_i^{-1}\cdots T_{j-2}^{-1}T_j^{-1} \cdots
T_{n-1}^{-1} \cdots T_i^{-1} & j<n \\[2mm]
t^{-n+i} x_n^{-1}\,I_{i,n-2}^{-1} & j=n \end{array} \right.
$$
Substituting this into (\ref{lapiz}) and hence into (\ref{boli})
along with (\ref{gafas}) (after using (\ref{reloj}) ) yields the
result. \hfill $\Box$

\begin{prop}
We have
$$
\left[ D_i, D_j \right] = 0 \hspace{3cm} 1\leq i, j \leq n
$$
\end{prop}
{\it Proof.}\quad Consider first the case $j=n$. In this case
\begin{eqnarray*}
\left[ D_i, D_n\right] &=& \left[ D_i, x_n^{-1}\left( 1-t^{n-1} Y_n\right)
\right] \\
&=& \left[D_i, x_n^{-1} \right] \left( 1-t^{n-1} Y_n\right) - t^{n-1}x_n^{-1}
\left[ D_i, Y_n \right] = 0
\end{eqnarray*}
thanks to (\ref{montreal}) and Lemma \ref{frio}. The general result now
follows using the representation (\ref{crep}) for $D_j$ in terms of $D_n$.
\hfill $\Box$

\setcounter{equation}{0}
\section{Raising/lowering operators}

The $q$-analogue of the raising operator $\Phi$ (recall (\ref{jphi}))
introduced by Knop and Sahi \cite{knop96c} is defined as
\begin{equation} \label{cansado.1}
\Phi_q := x_n\,T_{n-1}^{-1}\,\cdots\, T_2^{-1}\,T_1^{-1} =
t^{-n+i}\,T_{n-1}\,\cdots\,T_i\,x_i\,T_{i-1}^{-1}\,\cdots\,T_1^{-1}
\end{equation}
This operator enjoys the following properties

\begin{prop}
\begin{eqnarray*}
(a) \hspace{2cm}Y_j\,\Phi_q &=& \Phi_q\,Y_{j+1} \hspace{3cm} 1\leq j\leq n-1 \\
(b) \hspace{2cm}Y_n\, \Phi_q &=& q\,\Phi_q\,Y_1
\end{eqnarray*}
\end{prop}
{\it Proof.}\quad To prove (a), first note from (\ref{llave.1}) that
$$
Y_j\,\Phi_q = Y_j\,x_n\,T_{n-1}^{-1}\cdots T_1^{-1} = B_1 + B_2
$$
where
\begin{eqnarray*}
B_1 &=& x_n\, Y_j\, T_{n-1}^{-1}\cdots T_1^{-1}  \nonumber\\
&=& \Phi_q\,Y_{j+1} + (1-t^{-1})x_n\,T_{n-1}\cdots T_{j+1}^{-1}
T_{j-1}^{-1}\cdots T_1^{-1}\:Y_{j+1}
\end{eqnarray*}
and
\begin{eqnarray*}
B_2 &=& t^{-n+j}(1-t)x_n\,T_j\cdots T_{n-2}\,\omega\,T_1^{-1}\cdots
T_{j-1}^{-1}\:T_{n-1}^{-1}\cdots T_1^{-1} \nonumber\\
&=& t^{-n+j}(1-t)x_n\,T_{j-1}^{-1}\cdots T_1^{-1}\,\omega\,
T_1^{-1}\cdots T_j^{-1}
\end{eqnarray*}
A careful inspection of the second term occurring in $B_1$ shows that
it cancels with $B_2$, whence the result.

Similar considerations follow for (b), with the aid of (\ref{llave.2}).
\hfill $\Box$
\vspace{1mm}

The analogue of (\ref{ica.1}) is given by the following
\begin{cor}\label{phicor}
The operator $\Phi_q$ acts on non-symmetric Macdonald polynomials in
the following manner
$$
\Phi_q\; E_{\eta}(x;q,t) = t^{-\#\{\eta_i \leq \eta_1\} }\;
E_{\Phi\eta}(x;q,t)
$$
where $\Phi\eta := (\eta_2,\eta_3,\ldots,\eta_n,\eta_1+1)$.
\end{cor}
{\it Proof.}\quad From the previous Proposition, it is clear that
$\Phi_q\; E_{\eta}(x;q,t)$ is a contant multiple of 
$E_{\Phi\eta}(x;q,t)$ as they are both eigenfunctions of the operators
$Y_i$ with the same set of eigenvalues. The multiple is deduced by means
of examining the coefficient of the leading term $x^{\Phi\eta}$ in
the expansion of $\Phi_q\; E_{\eta}(x;q,t)$ with the aid of (\ref{t-inv})
and (\ref{t.action}). \hfill $\Box$
\vspace{1mm}

The definition of the lowering operator analogous to (\ref{jph}) makes 
use of the $q$-Dunkl operator introduced in Section 3: 
\begin{equation} \label{cansado.2}
\ph_q = T_1\,T_2\,\cdots\,T_{n-1}\,D_n = t^{n-i}\,T_1\,\cdots\,
T_{i-1}\,D_i\,T_i^{-1}\,\cdots\,T_{n-1}^{-1}
\end{equation}
This operator acts as a shift operator for the operators $Y_i$.

\begin{prop} \label{tarjeta}
\begin{eqnarray*}
(a) \hspace{2cm}Y_j\,\ph_q &=& \ph_q\,Y_{j-1} \hspace{3cm} 
2\leq j\leq n \\
(b) \hspace{2cm}Y_1\, \ph_q &=& q^{-1}\,\ph_q\,Y_n
\end{eqnarray*}
\end{prop}
{\it Proof.}\quad We begin with (a). Note that for $j>2$ we have
\begin{equation} \label{nevspeak}
T_j^{-1}\,T_1\cdots T_{n-1} = T_1\cdots T_{n-1}\, T_{j-1}^{-1} .
\end{equation}
Thus 
\begin{eqnarray}
Y_j\,\ph_q &=& t^{-n+j} T_j\cdots T_{n-1}\:\omega\: T_1^{-1}\cdots
T_{j-1}^{-1}\;(T_1\cdots T_n)\,D_n  \nonumber\\
&=& t^{j-2} D_1 \left( T_j \cdots T_{n-1}\:T_1^{-1}\cdots T_{n-2}^{-1}
\right)\: \omega\: T_1^{-1}\cdots T_{j-2}^{-1} \label{maquina}
\end{eqnarray}
where we have used (\ref{nevspeak}) along with Lemmas \ref{lemma.1}, 
\ref{lemma.2} and (\ref{tdefs.4}) to make the necessary
manipulations to get it into the above form. However we can rewrite the
term in the parenthesis occurring in (\ref{maquina}) as
$$
T_j \cdots T_{n-1}\:T_1^{-1}\cdots T_{n-2}^{-1} = 
T_1^{-1}\cdots T_{n-1}^{-1}\; T_{j-1}\cdots T_{n-1}
$$
which, when substituted back into (\ref{maquina}), yields 
the requisite result. The proof of (b) is almost immediate from
(\ref{horma.2}). \hfill $\Box$ 
\vspace{1mm}

As a consequence, we have the analogue of (\ref{ica.2})
\begin{cor} \label{phcor}
The action of $\ph_q$ on non-symmetric Macdonald polynomials is given
by
$$
\ph_q\, E_{\eta}(x;q,t) = t^{\#\{\eta_i<\eta_n\}}\;
\frac{d'_{\eta}(q,t)}{d'_{\ph\eta}(q,t)}\;E_{\ph\eta}(x;q,t)
$$
where $\ph\eta:= (\eta_n-1,\eta_1,\eta_2,\ldots,\eta_{n-1})$.
\end{cor}
{\it Proof.} \quad From Proposition \ref{tarjeta} we have that
$\ph_q\, E_{\eta}(x;q,t)$ is a multiple of $E_{\ph\eta}(x;q,t)$.
An examination of the leading term $x^{\ph\eta}$ in the expansion
of $\ph_q\, E_{\eta}(x;q,t)$ using (\ref{t.action}) and the explicit form
(\ref{explicit}) for $D_n$ tells us that
$$
\ph_q\, E_{\eta}(x;q,t) = t^{\#\{\eta_i<\eta_n\}}\;
(1-t^{n-1+\bar{\eta}_n} )\;E_{\ph\eta}(x;q,t) .
$$
However from Lemma \ref{props} we know that
$$
\frac{d'_{\eta}(q,t)}{d'_{\ph\eta}(q,t)} = 
1-t^{n-1+\bar{\eta}_n}
$$
whence the result.  \hfill $\Box$

\setcounter{equation}{0}
\section{Kernels}

We introduce the $q$-analogue of the Dunkl kernel ${\cal K}_A(x;y)$ by
\begin{equation} \label{ker.def}
{\cal K}_A(x;y;q,t) = \sum_{\eta} \frac{d_{\eta}(q,t)}{d'_{\eta}(q,t)
e_{\eta}(q,t)} \,E_{\eta}(x;q,t) E_{\eta}(y;q^{-1},t^{-1})
\end{equation}
This function reduces to ${\cal K}_A(x;y)$ as $q\rightarrow 1$ (although
${\cal K}_A(x;y;q,t)\neq {\cal K}_A(y;x;q,t)$ in general) and
satisfies generalizations of Theorem \ref{jackmain} and (\ref{belem}).
To establish these generalizations requires properties of the operators
$D_i$, $\Phi_q$ and $\ph_q$ obtained above, as well as some additional 
properties of the operators $T_i$ ($1\leq i\leq n-1$). The first such
property required relates to the action of the operators $T_i^{\pm 1}$
on the non-symmetric Macdonald polynomials \cite{noumi96c}
\begin{eqnarray}
T_i\,E_{\eta} &=& \left\{ \begin{array}{ll}
\left(\bfrac{t-1}{1-t^{-\de_{i\eta}}}\right)\,E_{\eta}
+ t\:E_{s_i\eta} & \eta_i < \eta_{i+1} \\
t\:E_{\eta} & \eta_i = \eta_{i+1} \\
\left(\bfrac{t-1}{1-t^{-\de_{i\eta}}}\right)\,E_{\eta} + 
\bfrac{(1-t^{\de_{i\eta}+1})(1-t^{\de_{i\eta}-1})}
{(1-t^{\de_{i\eta}})^2}\,E_{s_i\eta} & \eta_i > \eta_{i+1} 
\end{array}\right. \label{action1}\\
T^{-1}_i\,E_{\eta} &=& \left\{ \begin{array}{ll}
\left(\bfrac{t^{-1}-1}{1-t^{\de_{i\eta}}}\right)\,E_{\eta}
+ E_{s_i\eta} & \eta_i < \eta_{i+1} \\
t^{-1}\:E_{\eta} & \eta_i = \eta_{i+1} \\
\left(\bfrac{t^{-1}-1}{1-t^{\de_{i\eta}}}\right)\,E_{\eta} + 
t^{-1}\;\bfrac{(1-t^{\de_{i\eta}+1})(1-t^{\de_{i\eta}-1})}
{(1-t^{\de_{i\eta}})^2}\,E_{s_i\eta} & \eta_i > \eta_{i+1} 
\end{array}\right. \label{action2}
\end{eqnarray}
where $\delta_{i \eta} := \bar{\eta}_i - \bar{\eta}_{i+1}$.

Now define the involution \ $\inv{}$ \ as acting on operators or functions
by sending $q\rightarrow q^{-1}$, $t\rightarrow t^{-1}$.
The following lemma is the analogue of \cite[Lemma 3.7]{forr96d}
\begin{lemma} \label{symm}
Let $F(x,y)=\sum_{\eta} A_{\eta}\,E_{\eta}(x;q,t) E_{\eta}(y;q^{-1},t^{-1})$.
Then 
\begin{equation}
(T_i^{\pm 1})^{(x)}\,F(x,y) = 
\inv{T^{\mp 1}_i}^{\raisebox{-2mm}{\scriptsize$(y)$}} \,F(x,y) \label{cain1}
\end{equation}
if and only if the coefficients $A_{\eta}$ satisfy
\begin{equation} \label{llave}
A_{s_i\eta} = \left\{ \begin{array}{cc}
\bfrac{(1-t^{\de_{i\eta}})^2}{(1-t^{\de_{i\eta}+1})(1-t^{\de_{i\eta}-1})}
\; A_{\eta} & \eta_i > \eta_{i+1} \\[4mm]
\bfrac{(1-t^{\de_{i\eta}+1})(1-t^{\de_{i\eta}-1})}
{(1-t^{\de_{i\eta}})^2} \;A_{\eta} & \eta_i < \eta_{i+1} 
\end{array} \right.
\end{equation}
Moreover, these two conditions on $A_{\eta}$ are equivalent.
\end{lemma}
{\it Proof.}\quad Equation (\ref{cain1}) consists of two separate equations;
only $T_i^{(x)}\,F = \inv{T_i^{-1}}^{(y)}\,F$ will be established as the
other case follows from (\ref{t-inv}). The proof is similar to that given for
\cite[Lemma 3.7]{forr96d}. Split the sum in
$T_i^{(x)}\,F(x,y)$ according to whether $\eta_i < \eta_{i+1}$, 
$\eta_i = \eta_{i+1}$ or $\eta_i > \eta_{i+1}$. Apply (\ref{action1})
and collect coefficients of $E_{\eta}(x;q,t)$. Also, to work out
the action of $\inv{T^{-1}_i}^{\raisebox{-2mm}{\scriptsize$(y)$}}$ on
$E_{\eta}(y;q^{-1},t^{-1})$ (and hence on $F(x,y)$), set $t\rightarrow 
t^{-1}$ in (\ref{action2}). The two sides of (\ref{cain1}) are equal if 
and only if (\ref{llave}) holds. \hfill $\Box$

With this result at our disposal, the $q$-analogue of Theorem \ref{jackmain}
can now be given.

\begin{thm} \label{kernel}
The function $K_q(x;y)$ possesses the following properties:
\begin{eqnarray*}
(a)\qquad (T_i^{\pm 1})^{(x)}\,{\cal K}_A(x;y;q,t) &=& 
\inv{T_i^{\mp 1}}^{(y)}\,{\cal K}_A(x;y;q,t)\\
(b)\hspace{15mm} \ph_q^{(x)} \,{\cal K}_A(x;y;q,t) &=& 
\inv{\Phi_q}^{(y)}\,{\cal K}_A(x;y;q,t) \\
(c)\hspace{15mm} D_i^{(x)}\,{\cal K}_A(x;y;q,t) &=& y_i \,{\cal K}_A(x;y;q,t)
\end{eqnarray*}
\end{thm}
{\it Proof.}\\
(a) From Lemma \ref{props}, the constants $A_{\eta} = 
\bfrac{d_{\eta}(q,t)}{d'_{\eta}(q,t)e_{\eta}(q,t)}$ satisfy the
conditions of Lemma \ref{symm}, hence the result. \\
(b) From Lemma \ref{props} and Corollories \ref{phicor}, \ref{phcor}
we have
\begin{eqnarray*}
\ph_q^{(x)}\,{\cal K}_A(x;y;q,t) &=& \sum_{\eta} 
\frac{d_{\eta}(q,t)}{d'_{\eta}(q,t)e_{\eta}(q,t)}\,
t^{\#\{\eta_i<\eta_n\}}\frac{d'_{\eta}(q,t)}{d'_{\ph\eta}(q,t)}\,
E_{\ph\eta}(x;q,t) E_{\eta}(y;q^{-1},t^{-1}) \\
&=& \sum_{\nu} t^{\#\{\nu_i\leq \nu_1\}} 
\frac{d_{\Phi\nu}(q,t)}{d'_{\nu}(q,t)e_{\Phi\nu}(q,t)}\,
E_{\nu}(x;q,t) E_{\Phi\nu}(y;q^{-1},t^{-1}) \\
&=& \inv{\Phi_q}^{(y)}\,{\cal K}_A(x;y;q,t)
\end{eqnarray*}
(c) From (\ref{cansado.1}) we have
\begin{eqnarray*}
x_i &=& t^{n-i}\,T_i^{-1}\cdots T_{n-1}^{-1}\:\Phi_q\:
T_1\cdots T_{i-1} \\
&=& t^{-n+i}\,\inv{T_i^{-1}}\cdots \inv{T_{n-1}^{-1}}\:\inv{\Phi_q}\:
\inv{T_1}\cdots \inv{T_{i-1}}
\end{eqnarray*}
where the second form follows from applying the involution \ $\inv{}$
\ to the first form. Also, from (\ref{cansado.2}) 
$$
D_i = t^{-n+i}\,T_{i-1}^{-1}\cdots T_1^{-1}\:\ph_q\:T_{n-1}\cdots T_i
$$
The result now follows from these two expressions for $x_i$, $D_i$ 
by means of (a), (b) and the fact that operators acting on different
sets of variables commute. \hfill $\Box$
\vspace{.1cm}

It remains to present the analogue of (\ref{belem}). For the $q$-analogue
of Sym, Macdonald \cite{mac95} has introduced the operator
\begin{equation}
U^+ := \sum_{w \in S_n} T_w
\label{uplus}
\end{equation}
where $w = s_{i_1} \dots s_{i_p}$ $(1 \le i_1, \dots,i_p \le n-1)$ is the
reduced decomposition in terms of elementary transpositions of each element
of $S_n$ and
\begin{equation}
T_w = T_{i_1} \dots T_{i_p} 
\end{equation}
(the operators $T_i$ used by Macdonald satisfy $(T_i - t)(T_i + t^{-1})=0$
as distinct from (\ref{tdefs.1}); consequently in \cite{mac95} $U^+$ is defined
with $T_w$ multiplied by $t^{\ell(w)}$, where $\ell(w)$ is the length of
the permutation $w$, i.e. the number of elementary transpositions in its
reduced decomposition). As noted in \cite{mac95}, use of
(\ref{tdefs.1}) and (\ref{tdefs.2}) shows that
$
T_i U^+ = t U^+
$
which from the definition (\ref{pucallpa}) implies that $U^+f$ is symmetric 
in $x_1, \ldots, x_n$. In particular, for some proportionality constant
$a_\eta(q,t)$, we must have
\begin{equation}
U^+ E_\eta (x;q,t) = a_\eta(q,t) P_{\eta^+}(x;q,t)
\label{ue}
\end{equation}
where $P_{\eta^+}$ denotes the symmetric Macdonald polynomial normalized
so that the coefficient of the leading term is unity.

Our interest here is the action of $U^+$ on the kernel 
${\cal K}_A(x;y;q,t)$. For this we
require Theorem \ref{kernel} (a), (b) and (\ref{ue}) as well as the result of
the following lemma.

\begin{lemma} \label{dee}
Define
$$
(1-q) E_{0,m} := \sum_{i=m}^n A_{i,m} {(1 - \tau_i) \over x_i}
\quad {\rm where} \quad A_{i,m} := \prod_{j=m \atop j \ne i}^n
{tx_i - x_j \over x_i - x_j}.
$$
When acting on symmetric functions
\begin{equation}
\sum_{i=m}^n D_i = (1-q) E_{0,m}
\label{onep}
\end{equation}
\end{lemma}

\vspace{.2cm}
\noindent
{\it Proof.} \quad When acting on symmetric functions, we see from
(\ref{crep}) and (\ref{pucallpa}) that
\begin{equation}
D_j = T_j \dots T_{n-1} x_n^{-1} (1 - \tau_n).
\label{twoo}
\end{equation}
In particular $D_n = x_n^{-1}(1 - \tau_n)$ so (\ref{onep}) is true for
$m=n$. Thus by induction  (\ref{onep}) is equivalent to the statement that
\begin{equation}
D_{m-1} = (1 - q)(E_{0,m-1} - E_{0,m})
\label{two}
\end{equation}
Noting that
$$
A_{i,m-1} = \Big ( 1 + (t-1){x_i \over x_i - x_{m-1}} \Big )  A_{i,m}
$$
we see that (\ref{two}) can be rewritten to read
\begin{equation}
D_{m-1} = A_{m-1,m-1}{(1 - \tau_{m-1}) \over x_{m-1}} +
(t-1) \sum_{i=m}^n {x_i \over x_i - x_{m-1}} A_{i,m}
{(1 - \tau_i) \over x_i} := R_{m-1}
\label{twop}
\end{equation}
Since $R_n = D_n$, and from (\ref{twoo}) $T_j^{-1} D_j = D_{j+1}$, to
establish (\ref{twop}) it suffices to show
\begin{equation}
T_{m-1}^{-1} R_{m-1} = R_m
\label{three}
\end{equation}
This can be verified by direct calculation using (\ref{twop}) and
(\ref{t-inv}).
\hfill$\Box$

\vspace{.2cm}
We are now ready to calculate the action of $U^+$ on ${\cal K}_A(x;y;q,t)$.

\begin{prop}\label{propuf}
We have
\begin{equation}
U^{+(x)} {\cal K}_A(x;y;q,t) = [n]_t ! \: {}_0{\cal F}_0(x;y;q,t)
\label{tsymm}
\end{equation}
where $[n]_t ! := \prod_{i=1}^n (1-t^i)/(1-t)$, and 
with $\kappa$ denoting a partition and $b(\kappa) :=
\sum_{i=1}^n (i-1) \kappa_i$,
\begin{equation} \label{cano}
{}_0{\cal F}_0(x;y;q,t) := \sum_{\kappa} {t^{b(\kappa)} \over
d_\kappa'(q,t) P(1,t,\dots,t^{n-1};q,t)} P_\kappa(x;q,t) P_\kappa(y;q,t).
\end{equation}
\end{prop}

\noindent
{\it Proof. } \quad Applying $U^+$ to the $x$-variables in (\ref{ker.def})
gives
\begin{equation}
U^{+(x)} {\cal K}_A(x;y;q,t) = \sum_\eta {d_\eta(q,t) a_\eta(q,t) \over
d_\eta'(q,t) e_\eta(q,t)} P_{\eta^+}(x;q,t) E_\eta(y;q^{-1},t^{-1})
\label{prone}
\end{equation}
Repeating this operation and use of Theorem \ref{kernel} (a) shows
$$
U^{+(x)} {\cal K}_A(x;y;q,t) = \sum_{\eta^+} \alpha_{\eta^+}(q,t)  
P_{\eta^+}(x;q,t) P_{\eta^+}(y;q,t)
$$
where we have used the fact \cite{mac} that $ P_{\eta^+}(y;q^{-1},t^{-1})=
P_{\eta^+}(y;q,t)$. To specify $ \alpha_{\eta^+}$ sum  Theorem \ref{kernel}
(b) over $i$, apply $U^{+(x)}$ to both sides of the resulting equation and
commute its action to the right of $\sum_i D_i$ on 
the l.h.s.~(Lemma \ref{lemma.1} shows that this is
valid), and use Lemma \ref{dee} to substitute for $\sum_iD_i$ to show
$$
(1-q) E_{0,1}^{(x)}  \sum_{\eta^+} \alpha_{\eta^+}(q,t)  P_{\eta^+}(x;q,t)
P_{\eta^+}(y;q,t)
= p_1(y)  \sum_{\eta^+} \alpha_{\eta^+}(q,t)  P_{\eta^+}(x;q,t)
P_{\eta^+}(y;q,t)
$$
Recent results of Lassalle \cite[Theorems 3 and 5]{lass97a} give the action
of the operator $E_{0,1}^{(x)}$ on 
$ P_{\eta^+}(x;q,t)/P_{\eta^+}(1,\dots,t^{n-1};q,t)$
and an expression for the product $p_1(y) P_{\eta^+}(y;q,t)/d_{\eta}'(q,t)$
in terms of generalized binomial coefficients. These formulas imply
a recurrence for the quantity 
$\alpha_{\eta^+} d_{\eta^+}' P_{\eta^+}(1,\dots,t^{n-1};q,t)$
with the unique solution $\alpha_{\eta^+}d_{\eta^+}'
P_{\eta^+}(1,\dots,t^{n-1};q,t) = t^{b(\eta^+)}
\alpha_0$. The fact that $U^+ 1 = \sum_{w\in S_n} t^{\ell(w)} =[n]_t !$ 
gives $\alpha_0 = [n]_t !$.
\hfill$\Box$

\vspace{.2cm}
As an application of Proposition \ref{propuf} we can specify the 
proportionality constant $a_\eta$ in (\ref{ue}). For this we also 
require the formula \cite{mac95,noumi96c}
\begin{equation}
P_\kappa(y;q,t) = d_\kappa'(q,t) \sum_{\eta: \eta^+ = \kappa}
{1 \over d_\eta'(q,t)} E_\eta(y;q,t),
\label{rone}
\end{equation}
which is the analogue of (\ref{kirsten.1}).
We first substitute (\ref{prone}) on the l.h.s.~of (\ref{tsymm}), then
substitute for $P_{\eta^+}(y;q,t)$ on the l.h.s.~using 
$P_{\eta^+}(y;q,t) = P_{\eta^+}(y;q^{-1},t^{-1})={\rm RHS}\,(\ref{rone})
\Big |_{q \mapsto q^{-1}, t \mapsto t^{-1}}$. Next we note from the definition
(\ref{constants}) that 
$$
{d_{\eta^+}'(q^{-1},t^{-1}) \over d_{\eta}'(q^{-1},t^{-1}) } =
\prod_{s \in \eta^+, r \in \eta} t^{l(r) -l(s)} \;
{d_{\eta^+}'(q,t) \over d_{\eta}'(q,t) }
$$
and equate coefficients of $P_{\eta^+}(x;q,t)  E_\eta(y;q^{-1},t^{-1})$
on both sides to conclude
\begin{equation} \label{pavo}
a_\eta(q,t) = [n]_t !\; t^{\sum_{r\in \eta} l(r)} 
{e_\eta(q,t) \over P_{\eta^+}(1,\dots,t^{n-1};q,t) d_\eta(q,t)}.
\end{equation}
Substituting (\ref{pavo}) in (\ref{ue}), we obtain the $q$-analogue of
(\ref{kirsten.2}).

The function ${}_0{\cal F}_0(x;y;q,t)$ appears in an unpublished manuscript
of Macdonald \cite{macunp1}, as well as the recent work of Lassalle
\cite{lass97a}. Kaneko \cite{kaneko96a} introduced a similar function, with
$t^{b(\kappa)}$ in (\ref{cano}) replaced by $(-1)^{|\kappa|}q^{b(\kappa')}$.

\bibliographystyle{plain}

\end{document}